\documentclass[psamsfonts,reqno,11pt]{amsart}
\usepackage{graphicx}
\usepackage{graphics}
\usepackage{latexsym}
\usepackage{amsmath}
\usepackage{amssymb}
\usepackage[numbers,sort&compress]{natbib}
\usepackage[colorlinks, linkcolor=red,
anchorcolor= blue, citecolor=green]{hyperref}

\newtheorem{thm}{Theorem}[section]
\newtheorem{lem}[thm]{Lemma}
\newtheorem{eg}[thm]{Example}
\newtheorem{prop}[thm]{Proposition}
\newtheorem{cor}[thm]{Corollary}
\newtheorem{defn}[thm]{Definition}
\newtheorem{rem}[thm]{Remark}

\newtheorem{rem-eg}[thm]{Remark and Example}

\newtheorem{problem}[thm]{Problem}

\newcommand{\ti}{\tilde}

\newcommand{\OZ}{\mathbf{Z}}
\newcommand{\OX}{\mathbf{X}}
\newcommand{\OY}{\mathbf{Y}}
\newcommand{\OA}{\mathbf{A}}

\newcommand{\BH}{\mathbf{H}}
\newenvironment{prf}{{\noindent \textbf{Proof:} }}{\hfill $\Box$\medskip}
\def\sideremark#1{\ifvmode\leavevmode\fi\vadjust{\vbox
to0pt{\vss \hbox to 0pt{\hskip\hsize\hskip1em
\vbox{\hsize2cm\tiny\raggedright\pretolerance10000
\noindent#1\hfill}\hss}\vbox to8pt{\vfil}\vss}}}

\begin{document}

\title{On the numerical radius of Lipschitz operators in Banach spaces}
\author[Ruidong Wang]{Ruidong Wang}
\author[Xujian Huang]{Xujian Huang $^{*}$}
\author[Dongni Tan]{Dongni Tan}
\address[Ruidong Wang]{A Department of Mathematics, Tianjin University of Technology, Tianjin 300384, P.R. China}
\email{wangrdtjut@gmail.com}
\address[Xujian Huang]{A Department of Mathematics, Tianjin University of Technology, Tianjin 300384, P.R. China}
\email{huangxujian86@gmail.com}
\address[Dongni Tan]{A Department of Mathematics, Tianjin University of Technology, Tianjin 300384, P.R. China}
\email{tandongni0608@sina.cn}

\keywords{numerical radius, Lipschitz numerical index, Lipschitz operator.}

\subjclass[2010]{Primary 46B20; secondary 47A12.}

\thanks{$*$Corresponding author. \\ \indent  Second and third authors were supported by the National Natural Science Foundation
of China no. 11201339 and 11201338.}

\maketitle
\begin{abstract}
We study the numerical radius of Lipschitz operators on Banach spaces via the Lipschitz numerical index, which is an analogue of the numerical index in Banach space theory. We give a characterization of the numerical radius and obtain a necessary and sufficient condition for Banach spaces to have Lipschitz numerical index $1$. As an application, we show that real lush spaces and $C$-rich subspaces have Lipschitz numerical index $1$. Moreover, using the G$\hat{a}$teaux differentiability of Lipschitz operators, we characterize the Lipschitz numerical index of separable Banach spaces with the RNP. Finally, we prove that the Lipschitz numerical index has the stability properties for the $c_0$-, $l_1$-, and $l_\infty$-sums of spaces and vector-valued function spaces. From this, we show that the $C(K)$ spaces, $L_1(\mu)$-spaces and $L_\infty(\nu)$ spaces have Lipschitz numerical index $1$.
\end{abstract}

\section{Introduction}
\medskip

Our main goal in this paper is to study the numerical radius of Lipschitz operators by means of computing the Lipschitz numerical index of Banach spaces. The index is a constant related to the Lipschitz norm in the Banach algebra of all Lipschitz operators mapping a Banach space into the same space. Let us recall the relevant definitions. Let $\OX$ and $\OY$ be Banach spaces over the same coefficient field $\mathbb{R}$ or $\mathbb{C}$. A mapping $T$ from $\OX$ into $\OY$ is called {\it M-Lipschitz operator} if there exists a real constant $M> 0$ such that $$\|Tx-Ty\|\leq M\|x-y\|,\  \forall x,y \in \OX.$$ Let $Lip(\OX, \OY)$ denote the set of all Lipschitz operators from $\OX$ into $\OY$, and for
every $T \in Lip(\OX, \OY)$, let $\|T\|_L$ denote the minimum Lipschitz constant of $T$, i.e., $$\|T\|_L=\sup\{\frac{\|Tx-Ty\|}{\|x-y\|}: x, y\in \OX, x\neq y\}.$$ Then $\|\cdot\|_L$ is a semi-norm of the linear space $Lip(\OX, \OY)$. Let $Lip_0(\OX, \OY)$ denote the set of all Lipschitz operators from $\OX$ into $\OY$, which map $0$ to $0$. It is clear that $Lip_0(\OX, \OY)$ is a Banach space when it is equipped with the Lipschitz norm $\|\cdot\|_L$. Denote by $Lip_0(\OX)$ the space of all Lipschitz operators on $\OX$, which map $0$ to $0$. We shall sometimes refer to the Banach algebra $Lip_0(\OX)$ with the norm $\|\cdot\|_L$ as a Lipschitz operator algebra. For a real or complex Banach space $\OX$, we write $B_\OX, S_\OX$ and $\OX^*$ to denote its closed unit ball, its unit sphere and its dual space and denote the Banach algebra of all (bounded linear) operators on $\OX$ by $L(\OX)$. It is clear that $L(\OX)$ is a subalgebra of $Lip_0(\OX)$. For each $x\in \OX$, we define $$D(x)=\{x^*\in \OX^*: x^*(x)=\|x^*\|\cdot \|x\|=\|x\|^2\}.$$
For an operator $T\in Lip_0(\OX)$, its {\it numerical range} is defined as $$W(T)=\{\frac{(x-y)^*(Tx-Ty)}{\|x-y\|^2}: (x-y)^*\in D(x-y), \ x, y\in \OX, \ x\neq y\}.$$ and  its {\it numerical radius} is defined as $$\omega(T)=\sup\{|\lambda|: \lambda\in W(T)\}.$$ The {\it Lipschitz numerical index} of $\OX$ is the constant given by $$n_L(\OX)=\inf \{\omega(T): T\in Lip_0(\OX), \|T\|_L=1\}$$ or, equivalently, $$n_L(\OX)=\max \{k\geq 0: \|T\|\geq k\omega(T), T\in Lip_0(\OX)\}.$$
The definitions of Lipschitz numerical range is a generalization of the results of Zarantonello \cite{Z}, where only Hilbert spaces are considered.  This definition does work well for solving nonlinear functional equations and nonlinear partial differential equations (see \cite{D1}, \cite{Z}). Note that if $T\in L(\OX)$ is a bounded linear operator then $\omega(T)$ coincides with the usual numerical radius \cite{Wr}, i.e. $$\omega(T)=\sup\{|x^*(Tx)|: x^*\in D(x), x\in S_\OX\}.$$ The {\it numerical index} of $\OX$ is then given by  \begin{eqnarray*}n(\OX)=\inf \{\omega(T): T\in L(\OX), \|T\|=1\}.\end{eqnarray*} Obviously, $0\leq n_L(\OX)\leq n(\OX)\leq 1 $, $n_L(\OX) > 0$ means that the numerical radius is a norm on $Lip_0(\OX)$ equivalent to the Lipschitz norm and $n_L(\OX)=1$ if and only if the numerical radius and the Lipschitz norm coincide.
\medskip

Let us mention some facts concerning the numerical index which will be relevant to our discussion. The concept of numerical index of a Banach space was first suggested by G.Lumer in 1968. In the paper \cite{DMPW}, J. Duncan, C. McGregor, J. Pryce, and A. White determined the range of values of the numerical index, and showed that for complex (real) spaces the ranges of numerical index is the whole of the interval $[1/e, 1] ([0,1])$.  For every $T\in L(\OX)$, it is a well-known fact in the theory of numerical ranges (see \cite[\S 9]{BD} or \cite[Lemma 12]{L}) that $$\sup Re W(T)=\lim\limits_{t\rightarrow 0^+}\frac{\|I+tT\|-1}{t}$$
and so, $$\omega(T)=\max_{\alpha\in \mathbb{T}}\lim\limits_{t\rightarrow 0^+}\frac{\|I+t\alpha T\|-1}{t},$$ where $\mathbb{T}$ stands for the unit sphere of the base field $\mathbb{K}$ (=$\mathbb{R}$ or $\mathbb{C}$). Using the previous results, H. Bohnenblust and S. Karlin \cite{BK} or Glickfeld \cite{G} show that $n(\OX)\geq 1/e$ for any complex Banach space $\OX$. On the other hand, the numerical index of some classical Banach spaces have been calculated. For instance, it is known that if $\BH$ is a Hilbert space of dimension greater than $1$, $n(\BH) = 0$ in the real case and $n(\BH)=1/2$ in the complex case. The $L_1(\mu)$-spaces and their isometric preduals, including $C(K)$ spaces, function algebras and finite-codimensional subspaces of $C[0, 1]$ have numerical index $1$. For more information and background, we refer to the books by F.Bonsall and J.Duncan \cite{BD, BD1} and to the survey paper \cite{KMP} and references therein. Some known results about sums of Banach space and vector-valued function spaces appear in \cite{MP} and \cite{MV}. The numerical index of the $c_0$-, $l_1$-, or $l_\infty$-sum of a family of spaces coincides with the infimum of the numerical index of the spaces and the numerical indices of the vector-valued function spaces $C(K, \OX)$, $L_1(\mu,\OX)$ and $L_\infty(\nu, \OX)$ are equal to the numerical index of the Banach space $\OX$.
\medskip

In some sense, for a Banach space $\OX$ the Lipschitz operator algebra $Lip_0(\OX)$ is much bigger than the Banach algebra $L(\OX)$ of all bounded linear operators. Therefore, it may be natural to think the Lipschitz numerical index of $\OX$ should be less than or equal to the numerical index. Surprisingly, we show in this paper that for large classes of Banach spaces, the Lipschitz numerical index is equal to the numerical index.
\medskip

The outline of the paper is as follows.
\medskip

In Section 2 we generalize a well-known result of the numerical ranges of bounded linear operators to the Lipschitz operators and characterize
the numerical radius of Lipschitz operators. From this, we prove that the Lipschitz numerical index of a complex Banach space is greater than or equal to $1/e$.
We then give a necessary and sufficient condition for a Banach space to have Lipschitz numerical index $1$ and use this result to show that the real lush spaces also have Lipschitz numerical index $1$. Moreover, we establish that real and complex $C$-rich subspaces (specially, $c_0$, $C(K)$ spaces and finite-codimensional subspaces of $C[0, 1]$) have Lipschitz numerical index 1. As a consequence, all abelian $C^*$-algebras have Lipschitz numerical index 1.
\medskip

Section 3 is devoted to computing the Lipschitz numerical index through the use of the G$\hat{a}$teaux differentiability of Lipschitz operators. We give a characterization of the numerical radius of a Lipschitz operator, which is G$\hat{a}$teaux differentiable everywhere on a Banach space. Our main result is to show that the Lipschitz numerical index agrees with the numerical index on a separable Banach space with the Radon-Nikod$\acute{y}$m property. Using this result we show that for complex (real) spaces the ranges of Lipschitz numerical index is the whole of the interval $[1/e, 1] ([0,1])$.
\medskip

Finally, Section 4 is devoted to the study of the stability properties of the Lipschitz numerical index. We show that the Lipschitz numerical index of a $c_0$-, $l_1$-, and $l_\infty$-sums of a family of spaces is equal to the infimum Lipschitz numerical index of the summands. Thus the Lipschitz numerical index is stable for the $c_0$-, $l_1$-, and $l_\infty$-sums. For spaces of vector-valued functions we have the same results. Let $K$ be a compact Hausdorff space, $\mu$ a positive measure and $\nu$ a $\sigma$-finite measure. We prove that the spaces $C(K, \OX)$, $L_1(\mu,\OX)$ and $L_\infty(\nu, \OX)$ have the same Lipschitz numerical index as the Banach space $\OX$. From these results, a large family of classical spaces with Lipschitz numerical index $1$ are exhibited, namely $C(K)$ spaces, $L_1(\mu)$-spaces and $L_\infty(\nu)$-spaces.

\medskip

\section{Numerical radius of Lipschitz operators and Lush spaces}

Our first aim is to prove a formula connecting the Lipschitz norm with the numerical range of a Lipschitz operator. In fact, it is a generalization of the famous Lummer's lemma \cite[Lemma 12]{L} (see also \cite[\S 9]{BD}), which plays a fundamental role in the theory of the numerical index on Banach spaces. The proof is based on the one given in \cite[\S 9 Lemma 2]{BD} for the numerical range of bounded linear operators on a Banach space. We include it for the sake of completeness.

\begin{thm}\label{range}
Let $\OX$ be a Banach space. For each $T\in Lip_0(\OX)$
$$\sup Re W(T)=\lim\limits_{t\rightarrow 0^+}\frac{\|I+tT\|_L-1}{t}$$
and so $$\omega(T)=\max_{\alpha\in \mathbb{T}}\lim\limits_{t\rightarrow 0^+}\frac{\|I+t\alpha T\|_L-1}{t}.$$
\end{thm}
\begin{prf}
For each $x,y\in \OX$ with $x\neq y$ and $(x-y)^*\in D(x-y)$,
we have \begin{eqnarray*}\|I+tT\|_L&\geq& \frac{\|(I+tT)x-(I+tT)y\|}{\|x-y\|}\\ &\geq& Re\frac{(x-y)^*[(I+tT)x-(I+tT)y]}{\|x-y\|^2}
        \\ &=&1+Re\frac{t(x-y)^*(T x-T y)}{\|x-y\|^2}.\end{eqnarray*}
Thus $$\sup Re W(T)\leq\lim\limits_{t\rightarrow 0^+}\frac{\|I+tT\|_L-1}{t}.$$

On the other hand, let $\mu=\sup Re W(T)$. For each $x,y\in \OX$ with $x\neq y$ and $(x-y)^*\in D(x-y)$, we have
\begin{eqnarray*}\|(I-t T)x-(I-t T)y\| &\geq& Re \frac{(x-y)^*[(I-t T)x-(I-t T)y]}{\|x-y\|}
        \\ &=&\|x-y\|-t Re[\frac{(x-y)^*(T x-T y)}{\|x-y\|}]\\ &\geq& \|x-y\|(1-t\mu).\end{eqnarray*}
Hence, if we replace $x$ by $(I+t T)x$ and $y$ by $(I+t T)x$ then
\begin{eqnarray*}\|(I+t T)x-(I+t T)y\|&\leq& \frac{\|(I-t T)[(I+tT)x]-(I-t T)[(I+tT)y]\|}{1-t\mu}
        \\ &\leq& \frac{\|1-t^2T^2\|_L\|x-y\|}{1-t\mu}.\end{eqnarray*}
It follows that $\|I+tT\|_L\leq \|1-t^2T^2\|_L/(1-t\mu)$, and hence,
$$\lim\limits_{t\rightarrow 0^+}\frac{\|I+tT\|_L-1}{t}\leq \lim\limits_{t\rightarrow 0^+}\frac{\|I-t^2T^2\|_L/(1-t\mu)-1}{t}=\mu.$$
\end{prf}

The following variant of Theorem \ref{range} is very useful, which will be used later to get some stability properties of
the Lipschitz numerical index. We write $$\Pi=\{(x, y, f): x, y\in \OX, x\neq y,  f\in D(x-y)\}$$ and $\pi$ for the
projection from $\OX\times\OX\times\OX^*$ onto $\OX\times\OX$ defined by $\pi(x,y,f)=(x, y)$.
\begin{cor}\label{choose}
Let $Q\subset\Pi$ be such that $\pi(Q)$ is dense in $\OX\times\OX$. Then for each $T\in Lip_0(\OX)$,
$$\omega(T)=\sup \{\frac{|f(Tx-Ty)|}{\|x-y\|^2}: (x,y, f)\in Q \}.$$
Moreover, $\omega(T)$ may be determined by choosing one functional $f\in  D(x-y)$ at each point $(x,y)$ of a dense subset of $\OX\times\OX$.
\end{cor}

It is a celebrated result due to H.Bohnenblust and S.Karlin \cite{BK} (see also \cite{G}) that if $\OX$ is a complex Banach space, then $\|T\|\leq e\omega(T)$ for all  $T\in L(\OX)$. By Theorem \ref{range} we can generalize this result to the Lipschitz operators by showing that $\|T\|_L\leq e\omega(T)$ for all $T\in Lip_0(\OX)$.
Thus the numerical radius is always an equivalent Lipschitz norm in the complex case.

\begin{cor}\label{equivalent}
If $\OX$ is a complex Banach space, then $n_L(\OX)\geq 1/e$.
\end{cor}

We state other immediate consequences of Theorem \ref{range}, which are quite related to the so-called Daugavet equation
and alternative Daugavet equation. It is of interest to note the second statement, since it gives an equivalent condition which ensures that
the space has Lipschitz numerical index $1$.

\begin{cor}\label{numerical radius}
Let $\OX$ be a Banach space and $T\in Lip_0(\OX)$. Then

\smallskip\noindent

{\rm(a)} $\|I+T\|_L=1+\|T\|_L$ if and only if $\sup Re W(T)=\|T\|_L$.

\smallskip\noindent

{\rm(b)} $\max_{\alpha\in \mathbb{T}}\|I+\alpha T\|_L=1+\|T\|_L$ if and only if $\omega(T)=\|T\|_L$.
\medskip

Therefore, $\OX$ has Lipschitz numerical index $1$ if and only if all norm-one
operators $T$ in $Lip_0(\OX)$ satisfy $\max_{\alpha\in \mathbb{T}}\|I+\alpha T\|=2$.
\end{cor}

Next, we will study the Lipschitz numerical index of lush spaces by Corollary \ref{numerical radius}. The concept of lush space was introduced recently in \cite{BKMW}, which  has a geometrical property to ensure that the space has numerical index $1$. Some examples of lush spaces are $L_1(\mu)$-spaces and their isometric preduals, including $C(K)$ spaces, function algebras and finite-codimensional subspaces of $C[0, 1]$. We refer the reader to the papers \cite{BKMM, BKMW, KMMP} for more information and background on lush spaces.

\begin{defn}
A Banach space $\OX$ is said to be {\it lush} if for every $x,y \in S_\OX $ and every $\varepsilon> 0$, there is
a slice $$S=S(y^*, \varepsilon):=\{z\in B_\OX: Rey^*(z)>1-\varepsilon\}$$ with
$y^*\in S_{\OX^*}$ such that $y\in S$ and $dist(x, aco(S))<\varepsilon$, where $aco(S)$
denotes the absolutely convex hull of $S$.
\end{defn}

In order to simplify the writing, we introduce the following terminology. Let $\OX$ be a Banach space. For each $x,y\in \OX$, the {\it line segment} joining $x$ and $y$ is the set $$[x,y]=\{\lambda x+(1-\lambda) y: 0\leq \lambda\leq 1\}.$$ Let $A,B $ be subsets of $B_\OX$. We define the {\it join} of  $A$ and $B$ to be $$J(A, B)=\cup\{[x,y]: x\in A, y\in B \}$$ and we denote by  $J(A)$ the {\it join hull} of A which is the join of $A$ and itself. If $A$ is convex, then $J(A)$ is just $A$ itself.
\medskip

The following lemma is simple but very useful to compute the Lipschitz numerical index of Banach spaces.

\begin{lem}\label{join}
Let $\OX$ be a Banach space and $x,y \in \OX $ with $x\neq y$. If $\frac{x-y}{\|x-y\|}\in \overline{J(A)}$ for some set $A\subset B_\OX$, then for any $\varepsilon>0$ there exists $z\in \OX$ such that $$z-x\in \frac{\|x-y\|}{2}A \quad \mbox{and}  \quad \|y-z\|\leq (1+\varepsilon)\frac{\|x-y\|}{2}.$$
Moreover, if $\frac{x-y}{\|x-y\|}\in J(A)$, we can choose $z\in \OX$ such that $$z-x\in \frac{\|x-y\|}{2}A \quad \mbox{and}  \quad \|x-z\|=\|y-z\|=\frac{\|x-y\|}{2}.$$
\end{lem}

\begin{prf}
We write $z_0=\frac{x-y}{\|x-y\|}$ and choose $x_1,x_2\in A$,  $0\leq \lambda\leq 1$ such that $$\|z_0-(\lambda x_1+(1-\lambda)x_2)\|<\frac{\varepsilon}{2}.$$ We can assume that $\lambda \leq \frac{1}{2}$, then $z=x+\frac{\|x-y\|}{2}x_2$ is the desired element. Indeed, it is obvious that $$z-x=\frac{\|x-y\|}{2}x_2\in \frac{\|x-y\|}{2}A.$$ Moreover, set $v=\lambda x_1+(1-\lambda)x_2$. Since $\|2v-x_2\|\leq 1,$  it follows that
\begin{eqnarray*}\|z-y\|&=&\|x-y+\frac{\|x-y\|}{2}x_2\|=\frac{\|x-y\|}{2}\cdot\|2z_0-x_2\| \\&=& \frac{\|x-y\|}{2}\cdot\|(2z_0-2 v)+(2v-x_2)\|
 \leq (1+\varepsilon)\frac{\|x-y\|}{2}.\end{eqnarray*}
A similar argument proves the second part.
\end{prf}

We will prove that all real lush spaces have Lipschitz numerical index 1.

\begin{thm}\label{real lush}
Let $\OX$ be a real lush Banach space. Then $n_L(\OX)=1.$
\end{thm}

\begin{prf}
Let $T\in Lip_0(\OX)$ with $\|T\|_L=1$. From Corollary \ref{numerical radius}, it suffices to prove that $$\max\limits_{\alpha\in \mathbb{T}}\|I+\alpha T\|_L=2.$$ For every $\varepsilon>0$, there exist $x, y\in \OX$ such that $$\|T x-T y\|> (1-\varepsilon)\| x- y\|.$$
We apply the definition of lush spaces to $x_0=\frac{x-y}{\|x-y\|}$ and $y_0=\frac{Tx-Ty}{\|Tx-Ty\|}$
to obtain  $y_0^*\in S_{\OX^*}$ with $y_0\in S=S(y^*, \varepsilon)$ and $x_1\in S$,  $x_2\in -S$, $\lambda\geq 0$ such that $$\|x_0-(\lambda x_1+(1-\lambda) x_2)\|<\varepsilon.$$ The same proof of Lemma \ref{join} shows that there exists $z\in \OX$ such that $$z-x\in \frac{\|x-y\|}{2}(S\cup -S) \quad \mbox{and} \quad \|y-z\|\leq (1+\varepsilon)\frac{\|x-y\|}{2}.$$ Since moreover $y_0=\frac{Tx-Ty}{\|Tx-Ty\|}\in S=S(y^*, \varepsilon)$, we have $$|y^*(Tx-Ty)|>(1-\varepsilon)\|Tx-Ty\|>(1-2\varepsilon)\|x-y\|.$$
It follows that
\begin{eqnarray*}\max_{\alpha\in \mathbb{T}}\|I+\alpha T\|&\geq& \max_{\alpha\in\mathbb{T}}\frac{|y^*[(I+\alpha T)z-(I+\alpha T)x]|}{\|z-x\|} \geq \frac{|y^*(z-x)|+|y^*(T z-Tx)|}{\frac{\|x-y\|}{2}}\\ &\geq& 1-\varepsilon+\frac{|y^*(T x-T y)|-|y^*(T z-T y)|}{\frac{\|x-y\|}{2}}>2-6\varepsilon.
\end{eqnarray*}
Letting $\varepsilon \downarrow 0$, we deduce that $\max\limits_{\alpha\in \mathbb{T}}\|I+ \alpha T\|=2$ as desired.
\end{prf}

Let $K$ be a compact Hausdorff space in the real or complex case. We will give a direct proof that the $C(K)$ spaces
have Lipschitz numerical index 1 .

\begin{prop}\label{C(K)}
Let $K$ be a compact Hausdorff space. Then $n_L(C(K))=1.$
\end{prop}

\begin{prf}
Fix $T\in Lip_0(C(K))$ with $\|T\|_L=1$. Then for any $\varepsilon>0$ there exist $x, y \in C(K)$ such that $$\|Tx-Ty\|> (1-\varepsilon)\|x-y\|.$$
We may choose $s\in K$ with $x(s)-y(s)\neq 0$ such that \begin{eqnarray}|Tx(s)-T y(s)|>(1-\varepsilon)\|x-y\|.\end{eqnarray} Set $u=x-y$. Define $v\in C(K)$ and $g\in (C(K))^*$ given by $$v(t)=\frac{-u(t)}{\max\{|u(t)|, |u(s)|\}}, \ \forall \ t\in K \quad \mbox{and} \quad g(w)=\overline{-v(s)}w(s),\ \forall \ w\in C(K).$$ Then we clearly have $v\in S_{C(K)}$ and $g\in D(v)$. Let $z=x+\frac{\|x-y\|}{2}v$. Note that $$\omega(T)\geq \frac{|(\|z-x\|g)(Tz-Tx)|}{\|z-x\|^2}=\frac{|Tz(s)-Tx(s)|}{\|z-x\|}.$$ The key part is to show that $$|Tz(s)-Tx(s)|>(1-2\varepsilon)\|z-x\|.$$
Indeed,  for every $t\in K$,
$$z(t)-y(t)=x(t)-y(t)+\frac{\|x-y\|}{2}v(t)=u(t)(1-\frac{1}{2}\frac{\|u\|}{\max\{|u(t)|, |u(s)|\}});$$ If $|u(t)|>|u(s)|$, then
$$|z(t)-y(t)|=\Big||u(t)|-\frac{1}{2}\|u\|\Big|\leq \frac{1}{2}\|u\|=\frac{\|x-y\|}{2}.$$ If $|u(t)|\leq |u(s)|$, then $$|z(t)-y(t)|\leq\Big||u(s)|-\frac{1}{2}\|u\|\Big|\leq \frac{1}{2}\|u\|=\frac{\|x-y\|}{2}.$$ Thus $\|z-y\|\leq \frac{1}{2}\|x-y\|$,
and consequently, $$\|z-x\|=\|z-y\|=\frac{\|x-y\|}{2}.$$
Hence, from (1), we have
\begin{eqnarray*}|Tz(s)-Tx(s)|&\geq&|Tx(s)-Ty(s)|-|Tz(s)-Ty(s)|\\
&>&(1-\varepsilon)\|x-y\|-\|z-y\|=(1-2\varepsilon)\|z-x\|.\end{eqnarray*}
Hence,  $$\omega(T)\geq \frac{|Tz(s)-Tx(s)|}{\|z-x\|}>1-2\varepsilon.$$
Since $\varepsilon$ is arbitrary, it follows that $\omega(T)=1$. This completes the proof.
\end{prf}

We now present a wide class of subspaces of $C(K)$ which are introduced in the paper \cite{BKMW}, the so-called C-rich
subspaces. We will show that all C-rich subspaces of $C(K)$ have Lipschitz numerical index 1.

\begin{defn}\cite[Definition 2.3]{BKMW}
Let $K$ be a compact Hausdorff space. A closed subspace $X$ of
$C(K)$ is said to be \emph{C-rich} if for every nonempty open subset
$U$ of $K$ and every $\varepsilon>0$, there is a positive function
$h$ with $\|h\|=1$ and supp$(h)\subset U$ such that \mbox{dist}$(h,
X)<\varepsilon$.
\end{defn}
We next give some known examples of C-rich subspaces of $C(K)$.
\begin{eg}\label{eg}
\smallskip\noindent

{\rm(a)} Let $K$ be a compact Hausdorff space. Then $C(K)$ is C-rich.

\smallskip\noindent

{\rm(b)} If $K$ is a perfect compact space, then every finite-codimensional subspace of $C(K)$ is C-rich (see \cite[Proposition 2.5]{BKMW}).

\smallskip\noindent

{\rm(c)} If one considers $l_\infty$ as $C(\beta \mathbb{N})$, then every subspace of $l_\infty$ containing $c_0$ is C-rich.

\smallskip\noindent

{\rm(d)} Let $\Omega$ be locally compact, and let $\Omega_\infty=\Omega\cup \{\infty\}$ be the one-point compactification
of $\Omega$. Then $C_0(\Omega)$ is a C-rich subspace of $C(\Omega_\infty)$.
\end{eg}

\begin{defn}
We say that a Banach space $\OX$ has the {\it join-lush} property if for each $x,y \in S_\OX$ and $\varepsilon>0$, there exist $y^*\in S_{\OX^*}$ with $y \in S=S(y^*, \varepsilon)$, $x_1, x_2\in S$ and $\lambda\geq 0$ such that $$\|x-(\lambda \alpha_1 x_1+(1-\lambda) \alpha_2 x_2)\|<\varepsilon.$$ for some $\alpha_1, \alpha_2$ in $\mathbb{T}$.
\end{defn}

\begin{rem}\label{complex}
Note that the technique in the proof of Theorem \ref{real lush} is still valid in more general case. If $\OX$ is a Banach space with the join-lush property, then $\OX$ has Lipschitz numerical index 1.
\end{rem}

Due to the proof of \cite[Theorem 2.4]{BKMW}, every C-rich subspace of $C(K)$ is a lush space with the join-lush property. By the argument of Remark \ref{complex} we have the following result.

\begin{cor}\label{C-rich}
Let $K$ be a compact Hausdorff space and $\OX$ a C-rich subspace of $C(K)$. Then $n_L(\OX)=1$.
\end{cor}

By combing Corollary \ref{C-rich} with Example \ref{eg} (d) and Proposition \ref{C(K)}, the following is obtained.

\begin{cor}
If $\OA$ is an abelian $C^*$-algebra, then $n_L(\OA)=1$.
\end{cor}

\section{Numerical radius on separable Banach spaces with the RNP}

The G$\hat{a}$teaux differentiability of Lipschitz operators is a useful tool for studying the connection between Lipschitz operators and bounded linear operators. Let us recall some basic definitions. A mapping $T$ from an open set in a Banach space $\OX$ into a Banach space $\OY$ is said to be {\it G$\hat{a}$teaux differentiable} at a point $x_0$ if there is a bounded linear operator $S:\OX \rightarrow \OY$ such that for every $u\in \OX$, $$S(u)=\lim\limits_{t\rightarrow 0}\frac{T(x_0+tu)-T(x_0)}{t}.$$ The operator $S$ is called the {\it G$\hat{a}$teaux derivative} of $T$ at $x_0$ and denoted by $D_T(x_0)$. Clearly, if $T\in Lip_0(\OX, \OY)$ and $T$ is G$\hat{a}$teaux differentiable everywhere then $\|D_T(x)\|\leq \|T\|$ for every $x\in \OX$. A Banach space $\OY$ is said to have the {\it Radon-Nikod$\acute{y}$m property (RNP)} if every Lipschitz function $T : \mathbb{R} \rightarrow \OY$ is differentiable almost everywhere or equivalently every such $T$ has a point of differentiability. A Borel set $A$ in $\OX$ is said to be {\it Ar$\acute{o}$nszajn null} if for every sequence $\{x_n\}_n^\infty\subset \OX$ with a dense linear span, $A$ can be represented as a countable union $A=\cup A_n$ such that every line in the direction of $x_n$ meets $A_n$ is a set of measure $0$.
\medskip

We need the main existence theorem for G$\hat{a}$teaux derivatives of Lipschitz operators between Banach spaces.

\begin{thm}\label{null}
\rm{(\cite{A}, \cite{BL})} Let $\OX$ be a separable Banach space, and let $\OY$ be a Banach space with the RNP. Then every Lipschitz operator $T$ from an open set $U$
in $\OX$ into $\OY$ is G$\hat{a}$teaux differentiable outside an Ar$\acute{o}$nszajn null set.
\end{thm}

The following proposition will describe how the numerical radius of a G$\hat{a}$teaux differentiable Lipschitz operator relates to the one of its G$\hat{a}$teaux derivative.
\begin{prop}\label{radius}
Let $\OX$ be a Banach space, and let $T\in Lip_0(\OX)$ be G$\hat{a}$teaux differentiable everywhere. Then $$\omega(T)=\sup\{\omega(D_T(x)): x\in \OX\}.$$
\end{prop}

\begin{prf}
Let $x\in \OX$, $z\in S_\OX$ and $z^*\in D(z)$. By the definition of $D_T(x)$, we have
\begin{eqnarray*}
z^*(D_T(x)(z))=\lim\limits_{t\rightarrow 0}\frac{z^*(T(x+t z)-T(x))}{t}
=\lim\limits_{t\rightarrow 0}\frac{(t z)^*(T(x+tz)-T(x))}{t^2}.
\end{eqnarray*} It follows that $\omega(D_T(x))\leq \omega(T)$, and thus $$\sup\{\omega(D_T(x)): x\in \OX\}\leq \omega(T).$$

We shall prove the reverse inequality. For any $\varepsilon>0$, there exist $x_0,y_0\in \OX$ and $(x_0-y_0)^*\in D(x_0-y_0)$ such that
\begin{eqnarray}\frac{|(x_0-y_0)^*(Tx_0-Ty_0)|}{\|x_0-y_0\|^2}>\omega(T)-\frac{\varepsilon}{2}.\end{eqnarray} Set $z_0=\frac{x_0-y_0}{\|x_0-y_0\|}$, and fix $x\in [x_0, y_0]$. Since
$$\lim\limits_{t\rightarrow 0}\frac{T(x+t z_0)-T(x)}{t}=D_T(x)(z_0),$$ there exists sufficiently small $\eta_x>0$ such that $$\|D_T(x)(z_0)-\frac{T(x+tz_0)-T(x)}{t}\|<\frac{\varepsilon}{2}$$ for any $|t|<\eta_x$.
Let $V_x=(x-\eta_x/2z_0, x+\eta_x/2z_0)$ be such that $\{V_x: x\in[x_0, y_0]\}$ is an open cover of $[x_0, y_0]$. Choose finite points $\{x_i\}_{i=0}^n \subset [x_0, y_0]$ such that $x_n=y_0$ and $x_i\in (x_{i-1}, x_{i+1})$ for each $i=1,2, \cdots, n-1$.  So we get $$Tx_0-Ty_0=\sum\limits_{i=0}^{n-1}(Tx_i-Tx_{i+1}) \quad \mbox{and}  \quad \|x_0-y_0\|=\sum\limits_{i=0}^{n-1}\|x_i-x_{i+1}\|.$$ Let $z_0^*=\frac{(x_0-y_0)^*}{\|x_0-y_0\|}$. Then \begin{eqnarray}\frac{(x_0-y_0)^*(Tx_0-Ty_0)}{\|x_0-y_0\|^2}&=&\sum\limits_{i=0}^{n-1}\frac{(x_0-y_0)^*}{\|x_0-y_0\|}(\frac{Tx_i-Tx_{i+1}}{\|x_0-y_0\|})\nonumber
\\ &=&\sum\limits_{i=0}^{n-1}(\frac{\|x_i-x_{i+1}\|}{\|x_0-y_0\|})z_0^*(\frac{Tx_i-Tx_{i+1}}{\|x_i-x_{i+1}\|}).\end{eqnarray}
Hence, using (2) and (3), we may choose some $i_0\in \{0,1,\cdots, n-1 \}$ such that \begin{eqnarray}\frac{|z_0^*(Tx_{i_0}-Tx_{i_0+1})|}{\|x_{i_0}-x_{i_0+1}\|}>\omega(T)-\frac{\varepsilon}{2}.\end{eqnarray} We may assume that $\|x_{i_0}-x_{i_0+1}\|< \eta_{x_{i_0}}$ (otherwise, we can insert additional points between $x_{i_0}$ and $x_{i_0+1}$ to suit our purpose) and in view of the choice of $\eta_{x_{i_0}}$, we have \begin{eqnarray*}
\|\frac{T(x_{i_0})-T(x_{i_0+1})}{\|x_{i_0+1}-x_{i_0}\|}-D_T(x_{i_0})(z_0)\|\leq\frac{\varepsilon}{2}.
\end{eqnarray*}
This together with (4) proves that
\begin{eqnarray*}
|z_0^*\Big(D_T(x_{i_0})(z_0)\Big)| \geq \frac{|z_0^*(Tx_{i_0}-Tx_{i_0+1})|}{\|x_{i_0}-x_{i_0+1}\|}- \|\frac{T(x_{i_0})-T(x_{i_0+1})}{\|x_{i_0+1}-x_{i_0}\|}-D_T(x_{i_0})(z_0)\|\geq \omega(T)-\varepsilon.
\end{eqnarray*}
Since $\|z_0\|=1$, $z_0^*\in D(z_0)$ and $\varepsilon$ is arbitrary, we have $\sup\{\omega(D_T(x)): x\in \OX\}\geq \omega(T).$ We thus complete the proof.
\end{prf}

The proof of  Proposition \ref{radius} can be adapted to give other forms of the norm of Lipschitz operators.

\begin{prop}\label{norm}
Let $T$ be a Lipschitz operator from a Banach space $\OX$ to a Banach space $\OY$ such that $T(0)=0$.

\smallskip\noindent

{\rm(a)} Then $$\|T\|_L=\sup \{\frac{\|Tx-Ty\|}{\|x-y\|}: x\in \OX, y\in B(x,r), y\neq x\}$$ for every $r>0$, where $B(x,r)$ is the closed
ball with center $x$ and radius $r$ in $\OX$.

\smallskip\noindent

{\rm(b)} If $T$ is G$\hat{a}$teaux differentiable everywhere, then $$\|T\|_L=\sup\{\|D_T(x)\|: x\in \OX\}.$$
\end{prop}

\begin{prf}
(a) For every $x, y\in \OX$ with $x\neq y$, we select finite points $\{x_i\}_{i=0}^n \subset [x, y]$ with $x_0=x$ and $x_n=y$ such that
 $x_i\in (x_{i-1}, x_{i+1})$ for each $i=1,2, \cdots, n-1$. Then  \begin{eqnarray*}\frac{\|Tx-Ty\|}{\|x-y\|}=\sum\limits_{i=0}^{n-1}(\frac{\|x_i-x_{i+1}\|}{\|x-y\|})(\frac{Tx_i-Tx_{i+1}}{\|x_i-x_{i+1}\|}).\end{eqnarray*}
Since $\|x-y\|=\sum\limits_{i=0}^{n-1}\|x_i-x_{i+1}\|,$ by a simple convexity argument, there exists $i_0$ such that $$\frac{\|Tx_{i_0}-Tx_{i_0+1}\|}{\|x_{i_0}-x_{i_0+1}\|}\geq \frac{\|Tx-Ty\|}{\|x-y\|}.$$ The part (a) is then obtained by averaging the segment small enough.

(b) This part follows from the argument in Proposition \ref{radius} as well as the part (a).
\end{prf}

It is of interest to see whether or not the Lipschitz numerical index of a Banach space coincides with its numerical index. We will give some positive answers for separable Banach spaces with the RNP.

\begin{thm}\label{RNP}
Let $\OX$ be a separable Banach space satisfying the Radon-Nikod$\acute{y}$m property. Then $n_L(\OX)=n(\OX)$.
\end{thm}
\begin{prf}
We only need to prove that $n_L(\OX)\geq n(\OX)$. Let $T$ be in $Lip_0(\OX)$ with $\|T\|_L=1$.
We claim that for every $\varepsilon>0$, there exist $x_0, y_0\in \OX$ such that $T$ is G$\hat{a}$teaux differentiable almost everywhere in $[x_0,y_0]$ and $$\|T( x_0)-T( y_0)\|> (1-\varepsilon)\| x_0- y_0\|.$$ Indeed,  for the given $\varepsilon$, there exist $\ti x_0, \ti y_0\in X$ such that $$\|T(\ti x_0)-T(\ti y_0)\|> (1-\frac{1}{4}\varepsilon)\|\ti x_0-\ti y_0\|.$$
Let $\delta=\frac{\varepsilon}{8}\|\ti x_0-\ti y_0\|$. Then for each $x\in B(\ti x_0,\delta)$ and $y\in B(\ti y_0,\delta)$,
\begin{eqnarray*}
\|T(x)-T(y)\|&\geq&\|T(\ti x_0)-T(\ti y_0)\|-\|T(x)-T(\ti x_0)\|-\|T(y)-T(\ti y_0)\|\\ &>&(1-\frac{1}{2}\varepsilon)\|\ti x_0-\ti y_0\|
> \frac{1-\frac{1}{2}\varepsilon}{1+\frac{1}{2}\varepsilon}\|x-y\|>(1-\varepsilon)\|x-y\|.
\end{eqnarray*}
Moreover, the set $U=\{[x,y]: x\in B(\ti x_0,\delta), y\in B(\ti y_0,\delta)\}$ has nonempty interior in $\OX$.  By Theorem \ref{null}, there exist $x_0\in B(\ti x_0,\delta)$ and $y_0\in B(\ti y_0,\delta)$ such that $T$ is G$\hat{a}$teaux differentiable almost everywhere in $[x_0,y_0]$.

Let $z_0=\frac{y_0-x_0}{\|y_0-x_0\|}$ and $t_0=\|y_0-x_0\|$. Consider the mapping $f:[0,t_0]\rightarrow \OX$ given by
$f(t)=T(x_0+tz_0)$ for any $t\in [0,t_0]$. It is easily checked that $f$ is a Lipschitz function. Since $\OX$ has the RNP, by the proof of \cite[Theorem 5.21]{BL}, there exist a finite positive measure $\mu$ on $[0, t_0]$ which is absolutely continuous with respect to Lebesgue measure and a Bochner integrable $g:[0,t_0]\rightarrow \OX$ such that $$f(t)=\int_0^tg(s)d\mu(s)+f(0).$$
Then $f$ is G$\hat{a}$teaux differentiable and $D_f(t)=g(t)$ almost everywhere. On the other hand,
\begin{eqnarray*}
(1-\varepsilon)t_0\leq \|T(y_0)-T(x_0)\|=\|f(t_0)-f(0)\|\leq \int_0^{t_0}\|g(s)\|d\mu(s).
\end{eqnarray*}
So there exists a Borel set $A\subset [0,t_0]$ with $\mu(A)>0$ such that $\|g(s)\|\geq 1-\varepsilon$ for any $s\in A.$
Since $\mu(A)>0$ and $T$ is G$\hat{a}$teaux differentiable almost everywhere in $[x_0,y_0]$, it follows that there exists
$\ti t_0 \in A$ such that $T$ is G$\hat{a}$teaux differentiable at $x_0+\ti t_0z_0$. Hence
\begin{eqnarray*}
\|D_T(x_0+\ti t_0z_0)(z_0)\|&=&
\lim\limits_{s\rightarrow 0}\|\frac{T(x_0+\ti t_0z_0+s z_0)-T(x_0+\ti t_0z_0)}{s}\|\\
&=&\lim\limits_{s\rightarrow 0}\|\frac{f(\ti t_0+s)-f(\ti t_0)}{s}\|=\|D_f(\ti t_0)\|\\
&=&\|g(\ti t_0)\|\geq 1-\varepsilon.
\end{eqnarray*}
Taking $\hat{x}=x_0+\ti{t_0}z_0\in [x_0, y_0]$, we get an operator $D_T(\hat{x})$ in $L(\OX)$ satisfying $\|D_T(\hat{x})\|\geq 1-\varepsilon.$  Then the definition of $D_T(\hat{x})$ applies to show that  $$\omega(T)\geq \omega(D_T(\hat{x})) \geq (1-\varepsilon)n(\OX).$$
The desired inequality $n_L(\OX)\geq n(\OX)$ follows.
\end{prf}

Let $1\leq p<\infty$ and $(\Omega,\Sigma,\mu)$ be a $\sigma$-finite positive measure space. If $\Sigma$ is countably generated, then $L_p(\mu)$ is separable \cite[Prosition 3.4.5]{C}. Note that all reflexive spaces have the Radon-Nikod$\acute{y}$m property. Theorem \ref{RNP} gives the following consequence.

\begin{cor}\label{finite}

\smallskip\noindent

{\rm (a)} For any finite dimensional Banach space $\OX$, we have $n_L(\OX)=n(\OX)$.

\smallskip\noindent

{\rm (b)} Let $1<p<\infty$ and $(\Omega,\Sigma,\mu)$ be a $\sigma$-finite positive measure space. If $\Sigma$ is countably generated,
then $n_L(L_p(\mu))=n(L_p(\mu))$.

\smallskip\noindent

{\rm (c)} In the complex case, one has $n(l_2)=1/2$.
\end{cor}

From \cite[Theorems 3.5, 3.6]{DMPW}, J.Duncan, C.McGregor, J.Pryce, and A.White showed that for each $t\in [1/e,1]$ in the complex case (resp. $t\in [0, 1]$ in the real case) there is a two-dimensional complex (resp. real) space $\OX$ with $n(\OX)=t$. Then by Corollary \ref{finite} (a) and Corollary \ref{equivalent} we get the following result.
\begin{cor}
For a real Banach space $\OX$, $n_L(\OX)$ can be any
number in the interval $[0,1]$, while $\{n_L(\OX): \OX \ is \ a \ complex \ Banach \ space \} =[1/e, 1].$
\end{cor}

\section{ Stability properties of the Lipschitz numerical index.}
\medskip

The purpose of this section is to compute the Lipschitz numerical index of $c_0$-, $l_1$-, $l_\infty$-sums of Banach spaces and some vector-valued function spaces. It should be pointed out that the proof of the results are based on the ideas given in \cite{MP, MV}. However, we hope to convince the reader that the proof is not a straightforward adaptation of those given there, we need some techniques from the nonlinear theory.
\medskip

Given an arbitrary family $\{\OX_\lambda : \lambda\in \Lambda\}$ of Banach spaces, we denote by $[\oplus_{\lambda\in \Lambda}\OX_\lambda]_{c_0}$ (resp. $[\oplus_{\lambda\in \Lambda}\OX_\lambda]_{l_1}, [\oplus_{\lambda\in \Lambda}\OX_\lambda]_{l_\infty}$) the $c_0$-sum (resp. $l_1$-sum, $l_\infty$-sum) of the family. The sum of two spaces $\OX$ and $\OY$ is denoted by the simpler notation $\OX\oplus_\infty\OY$ or $\OX\oplus_1\OY$. For infinite countable sums of copies of a space $\OX$ we write $c_0(\OX), l_1(\OX), l_\infty(\OX)$.
\medskip

For our main results, we need two lemmas. The first one will be used repeatedly.

\begin{lem}\label{extension}
Let $\OX$ and $\OY$ be Banach spaces. If $x_1, x_2\in \OX$ and $y_1, y_2\in \OY$ such that $$\|y_1-y_2\|\leq M\|x_1-x_2\|$$ for some
$M>0$, then there
exists an $M$-Lipschitz operator $F\in Lip(\OX, \OY)$ such that $$F(x_1)=y_1\ \ and \ \ F(x_2)=y_2.$$
\end{lem}

\begin{prf}
The function $f_0: E_0=[x_1,x_2] \rightarrow \mathbb{R}$ defined by $f_0(x)=\|x-x_1\|$ is a Lipschitz function with $\|f_0\|_L=1$. By metric version of Tietze's extension theorem (see \cite[Theorem 1.5.6 (a)]{W}), there exists a function $f\in Lip(\OX, \mathbb{R}) $ such that $f|_{E_0}=f_0$ and $\|f\|_L=1$. Set $a=\|x_1-x_2\|$ and define $\phi:[0, a] \rightarrow \OY$ by $$\phi(t)=\frac{t}{a}y_2+(1-\frac{t}{a})y_1, \quad \forall \ t\in [0, a].$$ It is easy to see that $\phi$ is an  $M$-Lipschitz operator such that $$\phi(0)=y_1 \quad \mbox{and} \quad  \phi(a)=y_2.$$ Let $\pi_a: \mathbb{R}\rightarrow\mathbb{R}$ be a $1$-Lipschitz mapping defined by $$ \pi_a(z)= \left
\{ \begin{array}{ll}
|z|, & \mbox{ if }\; |z|\leq a;\\
a, & \mbox{ if }\;|z|>a.
\end{array}
\right.$$
Then $F=\phi\circ \pi_a \circ f: \OX \rightarrow \OY$ is the required mapping.
\end{prf}

\begin{lem}\label{infty}
Let $\OX,\OY$ and $\OZ$ be Banach spaces. Then for every $T\in Lip_0(\OX\oplus_\infty\OY, \OZ)$, there exist $u=(x_1,y_1), v=(x_2, y_2)\in \OX\oplus_\infty\OY$ with $\|u-v\|=\|x_1-x_2\|$ such that $$ \|Tu-Tv\|\geq(\|T\|_L-\varepsilon)\|u-v\|.$$
\end{lem}

\begin{prf}
We suppose that $T\in Lip_0(\OX\oplus_\infty\OY, \OZ)$ with $\|T\|_L=1$. Then for any $\varepsilon>0$ there exist $u, w\in \OX\oplus_\infty\OY$ such that $$\|Tu-Tw\|\geq(1-\frac{\varepsilon}{2})\|u-w\|.$$  It follows from Proposition \ref{norm} (a) that we can assume that $\|u-w\|\leq 1$. Since $B_{\OX\oplus_\infty\OY}$ is the join hull of $S_\OX \times B_\OY$, by Lemma \ref{join}, we may find $v\in \OX\oplus_\infty\OY$ such that  $$u-v\in S_\OX \times B_\OY \quad \mbox{and} \quad \|u-v\|=\|v-w\|=\frac{\|u-w\|}{2}.$$ A short computation shows that $\|T u-Tv\|\geq (1-\varepsilon)\|u-v\|$.
\end{prf}

\begin{prop}\label{sum}
Let $\{\OX_\lambda : \lambda\in \Lambda\}$ be a family of Banach spaces. If $\OZ$ is the $c_0$-, $l_1$-, or $l_\infty$-sum of the family, then $$n_L(\OZ)=\inf \{n_L(\OX_\lambda): \lambda\in \Lambda\}.$$
\end{prop}

\begin{prf}
It will be sufficient to prove the case that $\Lambda$ has just two elements, since in the general case given $\lambda_0\in \Lambda$, one clearly has $$\OZ=[\oplus_{\lambda \neq \lambda_0}\OX_\lambda]_{c_0}\oplus_\infty \OX_{\lambda_0} (resp. [\oplus_{\lambda \neq \lambda_0}\OX_\lambda]_{l_1}\oplus_1 \OX_{\lambda_0} \ or \ [\oplus_{\lambda \neq \lambda_0}\OX_\lambda]_{l_\infty}\oplus_\infty \OX_{\lambda_0}).$$

Let $\OZ$ denote either $\OX\oplus_\infty\OY$ or $\OX\oplus_1\OY$ for any Banach spaces $\OX$ and $\OY$.
We first check that $n_L(\OZ)\leq n_L(\OX)$. Let $S\in Lip_0(\OX)$ with $\|S\|_L=1$ and let $T\in Lip_0(\OZ)$ be given by $Tz=(Sx,0)$ for any $z=(x,y)\in \OZ$. Then $\|T\|_L = 1$, and given $\varepsilon>0$, we may find $z_1=(x_1,y_1), z_2=(x_2, y_2)\in \OZ$ with $x_1\neq x_2$ and $(f,g)\in D(z_1-z_2)$ such that \begin{eqnarray}\omega(T)-\varepsilon<\frac{|(f,g)(Tz_1-Tz_2)|}{\|z_1-z_2\|^2}=\frac{|f(Sx_1-Sx_2)|}{\|z_1-z_2\|^2}.\end{eqnarray}
Moreover, $(f,g)\in D(z_1-z_2)$ implies that \begin{eqnarray}f(x_1-x_2)+g(y_1-y_2)=\|f\|\|x_1-x_2\|+\|g\|\|y_1-y_2\|=\|z_1-z_2\|^2.\end{eqnarray}
Let $f_0=\frac{\|x_1-x_2\|}{\|f\|}f$. Then (5),(6) entail that $f_0\in D(x_1-x_2)$ and  $$\omega(T)-\varepsilon<\frac{|f(Sx_1-Sx_2)|}{\|z_1-z_2\|^2}\leq \frac{|f_0(Sx_1-Sx_2)|}{\|x_1-x_2\|^2}.$$  Consequently, we have $n_L(\OZ)-\varepsilon\leq \omega(T)-\varepsilon\leq\omega(S)$ and so $$n_L(\OZ)\leq n_L(\OX).$$ The same argument applies to $\OY$ giving
$n_L(\OZ)\leq \min\{n_L(\OX), n_L(\OY)\}$.

\medskip

Let $\OZ=\OX\oplus_\infty\OY$. To prove the reverse inequality, we consider the mapping $T\in Lip_0(\OZ)$ with $\|T\|_L=1$. Since $T$ can be rewritten as  $T=(T_1, T_2)$, where $T_1\in Lip_0(\OZ, \OX)$ and $T_2\in Lip_0(\OZ, \OY)$, we have $\|T\|_L=\max\{\|T_1\|_L, \|T_2\|_L\}$. We may assume that $\|T\|_L=\|T_1\|_L=1$. By Lemma \ref{infty}, there exist $z_1=(x_1,y_1), z_2=(x_2, y_2)\in \OZ$ with $\|x_1-x_2\|=\|z_1-z_2\|$ such that $$ \|T_1(z_1)-T_1(z_2)\|\geq(1-\varepsilon)\|z_1-z_2\|.$$  Lemma \ref{extension} gives a $1$-Lipschitz operator $F\in Lip(\OX, \OY)$ with $F(x_1)=y_1$ and $F(x_2)=y_2.$ Now consider the operator $S \in Lip_0(\OX)$ defined by $$S(u)=T_1\Big(u, F(u)\Big)-T_1\Big(0, F(0)\Big)\quad \forall \ u\in \OX .$$ It is clear that $ \|S(x_1)-S(x_2)\|\geq(1-\varepsilon)\|x_1-x_2\|$ and $\|S(u)-S(v)\|\leq\|u-v\|.$ This implies that $$1-\varepsilon\leq\|S\|\leq 1.$$ For all $x, \ti x \in \OX$ and $(x-\ti x)^*\in D(x-\ti x)$, we put $z=(x, F(x))$, $\ti z=(\ti x, F(\ti x))$ and $(z-\ti z)^*=((x-\ti x)^*,0)$. It is routine to check that $\|x-\ti x\|=\|z-\ti z\|$, $(z-\ti z)^*\in D(z-\ti z)$ and
$$(z-\ti z)^*(Tz-T \ti z)=(x-\ti x)^*(Sx-S\ti x).$$
It follows that $\omega(T)\geq \omega(S)\geq (1-\varepsilon) n_L(\OX)$, and hence $n_L(\OZ)\geq \min\{n_L(\OX), n_L(\OY)\}$.
\medskip

Let $\OZ=\OX\oplus_1\OY$. We shall prove the reverse inequality $n_L(\OZ)\geq \min\{n_L(\OX), n_L(\OY)\}$. Fix $T\in Lip_0(\OZ)$ with $\|T\|_L=1$. Given $0<\varepsilon<1$, there exist two elements $z_1=(x_1, y_1), z_2=(x_2, y_2)\in \OZ$ such that $$\|T(x_1, y_1)-T(x_2, y_2)\|\geq(1-\varepsilon)(\|x_1-x_2\|+\|y_1-y_2\|).$$ Thus we have $$\|T(x_1, y_1)-T(x_1, y_2)\|\geq(1-\varepsilon)\|y_1-y_2\| \quad \mbox{or} \quad \|T(x_1, y_2)-T(x_2, y_2)\|\geq (1-\varepsilon) \|x_1-x_2\|.$$ Since the argument works whenever one of the inequalities is true, we may assume that the second inequality holds.  Then define $T_{y_2}:\OX\rightarrow\OX\oplus_1\OY$ by
$$T_{y_2}(x)=T(x, y_2)-T(0, y_2), \quad \forall x\in \OX,$$ and write $T_{y_2}=(A, B)$, where $A\in Lip_0(\OX)$ and $B\in Lip_0(\OX, \OY)$. Let $x_0\in S_\OX$ and $y^*\in S_{\OY^*}$ be such that $$Ax_1-Ax_2=\|Ax_1-Ax_2\|x_0 \quad \mbox{and} \quad y^*(Bx_1-Bx_2)=\|Bx_1-Bx_2\|,$$ and define an operator $S\in Lip_0(\OX)$ by $$S(x)=A(x)+y^*(Bx)x_0, \quad \forall x\in \OX.$$ Then it is easily checked that $\|S(x_1)-S(x_2)\|\geq(1-\varepsilon)\|x_1-x_2\|$, and thus $$\|S\|_L\geq 1-\varepsilon.$$ Finally,  for all $\ti x_1, \ti x_2\in \OX$ and $(\ti x_1- \ti x_2)^*\in D(\ti x_1- \ti x_2)$, we set $\ti z_1=(\ti x_1, y_2), \ti z_2=(\ti x_2, y_2) \in \OZ$ and $f=\Big((\ti x_1- \ti x_2)^*, (\ti x_1- \ti x_2)^*(x_0)y^*\Big)\in \OZ^*$. It is easy to see that
$\|\ti z_1-\ti z_2\|=\|\ti x_1- \ti x_2\|$, $f\in D(\ti z_1- \ti z_2)$ and
$$f(T \ti z_1- T \ti z_2 )=f\Big(T_{y_2}(\ti x_1)-T_{y_2}(\ti x_2)\Big)=(\ti x_1- \ti x_2)^* ( S \ti x_1-S \ti x_2).$$ It follows that $$\omega(T)\geq \omega(S) \geq (1-\varepsilon) n_L(\OX).$$ The desired inequality $n_L(\OZ)\geq \min\{n_L(\OX), n_L(\OY)\}$ follows.
\end{prf}

As an application of Proposition \ref{sum}, we show that the Lipschitz numerical index is stable under $c_0$-, $l_1$-, and
$l_\infty$-sums.

\begin{cor}
For every Banach space $\OX$, $$n_L(c_0(\OX))=n_L(l_1(\OX))=n_L(l_\infty(\OX))=n_L(\OX).$$
In particular, $n_L(c_0)=n_L(l_1)=n_L(l_\infty)=1$ in the real or complex case.
\end{cor}

Next, we will discuss the stability property of spaces of vector-valued functions.
Let us recall some notation. Given a compact Hausdorff space $K$ and a Banach space $\OX$, we denote by $C(K,\OX)$ the Banach space of all continuous functions from $K$ into $\OX$, endowed with the supremum norm. If $(\Omega, \Sigma, \mu)$ is a positive measure space, $L_1(\mu, \OX)$ is the Banach space of all Bochner-integrable functions $f: \Omega\rightarrow \OX$ with $$\|f\|_1=\int_\Omega \|f(t)\| d\mu(t).$$ If $(\Omega, \Sigma, \mu)$ is $\sigma$-finite, then $L_\infty(\mu, \OX)$ stands for the space of all Bochner integrable functions $f$ from $\Omega$ into a Banach space $\OX$, endowed with its natural norm $$\|f\|_\infty= \inf\{\lambda\geq 0: \|f(t)\|\leq \lambda \ a.e. \} .$$ We refer to \cite{D} for more background information.

We shall generalize the fact $n_L(c_0(\OX))=n_L(\OX)$ to the space of vector-valued continuous functions.

\begin{thm}\label{compact}
Let $K$ be a compact Hausdorff space and $\OX$ a Banach space. Then, $$n_L(C(K,\OX))=n_L(\OX).$$
\end{thm}

\begin{prf}
We first show that $n_L(C(K,\OX)) \geq n_L(\OX)$. Let $T$ be in $Lip_0(C(K,\OX))$ with $\|T\|_L=1$, and the procedure is to prove that $\omega(T)\geq n_L(\OX)$.
For each $t\in K$, we define $T_t \in Lip_0(C(K,\OX), \OX)$ by $T_t(h)=(Th)(t)$ for any $h\in C(K,\OX)$. Then $$\|T\|_L=\sup\{\|T_t\|_L:t\in K\}.$$
Given $\varepsilon>0$, we may find $t_0\in K$ such that $\|T_{t_0}\|_L>1-\varepsilon$. Thus there exist $f,g\in C(K,\OX)$ such that $$\|T_{t_0}f-T_{t_0}g\|>(1-\varepsilon)\|f-g\|.$$ Set $f_0=\frac{f-g}{\|f-g\|}$, $C(t_0)=\{\varphi\in C(K): \varphi (K) \subset [0,1], \varphi(t_0)=1\}$ and $$\mathcal{A}(f_0)=\{(1-\varphi) f_0+\varphi x: x\in S_\OX, \varphi\in C(t_0)\}.$$ By the proof of \cite[Theorem 5]{MP}, we have $f_0 \in \overline{J(\mathcal{A}(f_0))}$ (Indeed, given $\varepsilon>0$, we write $z_0=f_0(t_0)$ and $z_0 =\lambda x_1+(1-\lambda)x_2$ with $0\leq \lambda\leq 1$, $x_1,x_2 \in S_\OX$. Then there exists a continuous function $\varphi : K \rightarrow[0,1]$ such that $\varphi(t_0)=1$ and $\varphi(t)=0$ if $\|f_0(t)-z_0\|\geq \varepsilon$ and consider the functions $$f_j =(1-\varphi) f_0+\varphi x_j \in \mathcal{A}(f_0) \ \ (j= 1,2). $$ Then $\|f_0-(\lambda f_1+ (1-\lambda) f_2)\|<\varepsilon$).  We see from Lemma \ref{join} that, there is a function $h\in C(K,\OX)$ such that $$ h-f \in \frac{\|f-g\|}{2}\mathcal{A}(f_0) \quad  \mbox{and} \quad \|h-g\|\leq (1+\varepsilon)\frac{\|f-g\|}{2}.$$
Hence, we have $$\|T_{t_0}h-T_{t_0}f\|\geq \|T_{t_0}f-T_{t_0}g\|-\|h-g\|\geq (1-3\varepsilon)\|h-f\|.$$
Now write $x_0=h(t_0)$, $y_0=f(t_0)$ and so $\|x_0-y_0\|=\|h-f\|.$ By Lemma \ref{extension}, there exists a $1$-Lipschitz operator $F\in Lip(\OX, C(K, \OX))$ such that $$F(x_0)=h\ \quad \mbox{and} \quad \ F(y_0)=f.$$ Next we find a new continuous function $\phi : K \rightarrow[0,1]$ such that $\phi(t_0)=1$ and $\phi(t)=0$ if
$$\|h(t)-x_0\|\geq \varepsilon\|h-f\| \quad \mbox{and} \quad \|f(t)-y_0\|\geq \varepsilon\|h-f\|,$$ and denote
$$\Phi(x)=(1-\phi)F(x)+\phi x\in C(K, \OX), (x\in \OX).$$ It follows that
\begin{eqnarray*}\|T_{t_0}(\Phi(x_0))-T_{t_0}(\Phi(y_0))\|&=&\|T_{t_0}((1-\phi)h+\phi x_0)-T_{t_0}((1-\phi)f+\phi y_0)\|
\\ &\geq& \|T_{t_0}h-T_{t_0}f\|-\|\phi h-\phi x_0\|-\|\phi f-\phi y_0\|
\\ &\geq& (1-3\varepsilon)\|h-f\|-\varepsilon\|h-f\|-\varepsilon\|h-f\|
\\ &=&(1-5\varepsilon)\|h-f\|=(1-5\varepsilon)\|x_0-y_0\|.\end{eqnarray*}
We consider the operator $S \in Lip_0(\OX)$ given by $$S(x)=T_{t_0}(\Phi(x))-T_{t_0}(\Phi(0))$$ for any $x\in \OX$. We can easily check that
$\|Sx_0-Sy_0\|\geq (1-5\varepsilon)\|x_0-y_0\|$ and $$\|Sx-Sy\|\leq \|\Phi(x)-\Phi(y)\|=\|x-y\|$$ for all $x, y \in \OX.$
Hence $1-5\varepsilon\leq\|S\|\leq 1$. For each $x, y\in \OX$ and $(x-y)^*\in D(x-y)$, we have $(x-y)^*\circ\delta_{t_0} \in D(\Phi(x)-\Phi(y))$
and $$(x-y)^*\circ\delta_{t_0}(T\Phi(x)-T\Phi(y))=(x-y)^*(T_{t_0}\Phi(x)-T_{t_0}\Phi(y))=(x-y)^*(Sx-Sy).$$ Hence,
$$\omega(T)\geq \omega(S)\geq (1-5\varepsilon)n_L(\OX)$$ and thus $n_L(C(K, \OX))\geq n_L(\OX)$.

To prove the reverse inequality, for every $S \in Lip_0(\OX)$ with $\|S\|_L=1$ we define $T\in Lip_0(C(K,\OX))$ given by
$$(Tf)(t)=S(f(t)) \quad \,  t\in K, f\in C(K, \OX).$$ Then $\|T\|_L=1$ and $\omega(T)\geq n_L(C(K,\OX))$. Let $\OZ=C(K, \OX)$. For every $t\in K$, we set $$\mathcal{A}_t=\{(f,g): \|f(t)-g(t)\|=\|f-g\|\} \subset \OZ \times \OZ$$ and $$Q=\{(f, g, x^*\circ \delta_t): (f,g)\in \mathcal{A}_t, x^*\in D(f(t)-g(t)), t\in K \}\subset \OZ\times \OZ\times \OZ^*.$$ Then the mapping $\pi: \OZ\times\OZ\times\OZ^* \rightarrow \OZ\times\OZ$ given by $\pi(z_1, z_2, z^*)=(z_1, z_2)$ satisfies $\pi(Q)=\OZ\times\OZ$. By Corollary \ref{choose} the numerical radius of $T$ is given by \begin{eqnarray*}\omega(T)=\sup\{\frac{|(x^*\circ \delta_t)(T f-T g)|}{\|f-g\|^2}: (f, g, x^*\circ \delta_t) \in Q \}.\end{eqnarray*}  Therefore, given $\varepsilon>0$, we can find $(f, g, x^*\circ \delta_t)\in Q$ such that $$\omega(T)-\varepsilon< \frac{|x^*[(Tf)(t)-(Tg)(t)]|}{\|f-g\|^2}=\frac{|x^*[S(f(t))-S(g(t))]|}{\|f(t)-g(t)\|^2}.$$
It follows that $\omega(S)\geq \omega(T)-\varepsilon\geq n_L(C(K, \OX))-\varepsilon$ and so $n_L(\OX)\geq n_L(C(K, \OX))$.
\end{prf}

We will generalize the fact that $n_L(l_1(\OX))= n_L(\OX)$ as follows.
\begin{thm}\label{positive measure}
Let $(\Omega, \Sigma,\mu)$ be a positive measure space,  and let $\OX$ be a Banach space. Then $$n_L(L_1(\mu, \OX))=n_L(\OX).$$
\end{thm}

The proof will depend on the following two lemmas.
\begin{lem}\label{finite range}
Let $\OX$ be a Banach space, and let $S$ be in $Lip_0(\OX)$. If  $T\in Lip_0(l_1^n(\OX))$ defined by $T(x)=(Sx_1, Sx_2, \cdots, Sx_n)$ for each $x=(x_1, x_2, \cdots, x_n) \in l_1^n(\OX)$, then $$\omega(T)=\omega(S).$$
\end{lem}

\begin{prf}
It is enough to prove that $\omega(T)\leq\omega(S)$ since the converse is obvious.
For every $\varepsilon>0$, we may find $x_0=(\ti x_1, \ti x_2, \cdots, \ti x_n)$,  $ y_0=(\ti y_1, \ti y_2, \cdots, \ti y_n)$ in $l_1^n(\OX)$ with $\ti x_i\neq \ti y_i$, for all $i=1,2, \cdots, n$ and $ z^*=(z_1^*, z_2^*, \cdots, z_n^*) \in D(x_0-y_0)$ such that
 \begin{eqnarray}\frac{ |z^*(T  x_0-T y_0)|}{\| x_0- y_0\|^2}\geq \omega(T)-\varepsilon.\end{eqnarray}
It follows from $z^* \in D(x_0-y_0)$ that $$z_i^*(\ti x_i-\ti y_i)=\|z_i^*\|\|\ti x_i-\ti y_i\| \quad
\mbox {and} \quad \|z_i^*\|=\|x_0-y_0\|$$ for all $i=1,2, \cdots, n.$
This implies that each $\ti z_i^*=\frac{\|\ti x_i -\ti y_i\|}{\|x_0-y_0\|}z_i^*\in D(\ti x_i-\ti y_i)$ and \begin{eqnarray} \frac{ z^*(T  x_0-T  y_0)}{\| x_0- y_0\|^2}&=&\sum\limits_{i=1}^n\frac{ z_i^*(S \ti x_i-S\ti y_i)}{\| x_0- y_0\|^2} =\sum\limits_{i=1}^{n}(\frac{\|\ti x_{i}- \ti y_{i}\|}{\|x_0-y_0\|})\frac{\ti z_i^*(S \ti x_{i}- S\ti y_{i})}{\|\ti x_{i}- \ti y_{i}\|^2}.\end{eqnarray}
By (7) and (8), there exists $i_0\in \{1, 2, \cdots, n\}$ such that $$\frac{\ti z_{i_0}^*(S \ti x_{i_0}- S\ti y_{i_0})}{\|\ti x_{i_0}- \ti y_{i_0}\|^2}\geq \omega(T)-\varepsilon.$$ Therefore,
$\omega(T)\leq\omega(S).$
\end{prf}

Let us recall some definitions. Let $(\Omega, \Sigma,\mu)$ be a positive finite measure space and $\OX$ a Banach space.
For each $A\in \Sigma$ with $0<\mu(A)<\infty$ we define a contractive linear operator $\Theta_A: L_1(\mu, \OX)\rightarrow \OX$ by $$\Theta_A(f)=\frac{1}{\mu(A)}\int_A f d\mu \quad \forall \ f\in L_1(\mu, \OX).$$ With this notation, for each partition $\pi$ of $\Omega$
(into a finite set of disjoint members of $\Sigma$) we define a contractive linear
operator $E_\pi:L_1(\mu, \OX)\rightarrow L_1(\mu, \OX)$ given by $$E_\pi(f)=\sum\limits_{A\in\pi}\Theta_A(f)\chi_A, \quad \forall f\in L_1(\mu, \OX).$$

\begin{lem}\cite[Lemma III.2.1]{D}\label{partition}
Let $(\Omega, \Sigma,\mu)$ be a positive finite measure space and $\OX$ a Banach space. If the partitions are directed by refinement,
then $$\lim\limits_\pi\|E_\pi(f)-f\|_1=0, \quad \forall \ f\in  L_1(\mu, \OX).$$
\end{lem}
\medskip
{\bf Proof of Theorem \ref{positive measure}.}
We adapt the proof of \cite[Theorem 8]{MP} to prove this result.
According to the fact that $L_1(\mu, \OX)$ is isometrically isomorphic to an $l_1$-sum of spaces $L_1(\mu_i, \OX)$ for suitable finite measures $\mu_i$, and together with Proposition \ref{sum}, we can suppose that $\mu$ is finite measure.

For each partition $\pi$ of $\Omega$, we denote by  $\OY_\pi$ the range of $E_\pi$. Since $\OY_\pi$ is isometric to a finite $l_1$-sum of copies of $\OX$, by Proposition \ref{sum} we have $n_L(\OY_\pi)=n_L(\OX)$.

We shall prove the inequality $n_L(L_1(\mu, \OX))\geq n_L(\OX)$. Let $T\in Lip_0(\OX)$ with $\|T\|_L=1$ and $\varepsilon>0$. By the previous consideration it suffices to find a partition $\pi_0$ of $\Omega$ and a mapping $T_{\pi_0}\in Lip_0(\OY_{\pi_0})$ such that $$\|T_{\pi_0}\|_L\geq 1-\varepsilon \quad \mbox{and} \quad \omega(T)\geq \omega(T_{\pi_0})\geq (1-\varepsilon) n_L(\OY_{\pi_0})=(1-\varepsilon)n_L(\OX).$$ For the given $T$ and $\varepsilon>0$, we may find $f,g\in L_1(\mu, \OX)$ such that $$\|Tf-Tg\|>(1-\frac{1}{4}\varepsilon)\|f-g\|.$$
Using Lemma \ref{partition}, we obtain a partition $\pi_0$ of $\Omega$ into a finite
family of disjoint measurable sets with positive measure satisfying $$\max\{\|f-E_{\pi_0} f\|, \ \|Tf-E_{\pi_0} Tf\|, \ \|g-E_{\pi_0} g\|, \ \|Tg-E_{\pi_0} Tg\|\}\leq \frac{1}{4}\varepsilon\|f-g\|.$$  Let $T_{\pi_0}=E_{\pi_0} T E_{\pi_0} \in Lip_0(\OY_{\pi_0})$ be the restriction to $\OY_{\pi_0}$ of the operator $E_{\pi_0} T$. Then it is easily checked that $\|T_{\pi_0}\|_L\geq 1-\varepsilon$. Moreover, for each $f, g\in \OY_{\pi_0}$ and $(f-g)^*\in D(f-g)$, we have $h^*=E_{\pi_0}^* (f-g)^*\in D(f-g)$ and $$h^*(Tf-Tg)=(f-g)^*(T_{\pi_0}f-T_{\pi_0}g),$$ which shows that $\omega(T)\geq \omega(T_{\pi_0})$ as required.
\medskip

Next, we check the reverse inequality $n_L(L_1(\mu, \OX))\leq n_L(\OX)$. For every $S \in Lip_0(\OX)$ with $\|S\|_L=1$, we define $T\in Lip_0(L_1(\mu, \OX))$ by $(Tf)(t)=S(f(t))$ for each $t\in \Omega, \ f\in L_1(\mu, \OX).$ Then $\|T\|_L=1$ and $\omega(T)\geq n_L(L_1(\mu, \OX))$. Given $\varepsilon>0$, by Corollary \ref{choose} and Lemma \ref{partition} we can choose a finite partition $\pi$ of $\Omega$ such that $$\|T_\pi\|_L\geq 1-\varepsilon \quad \mbox{and} \quad \omega(T_\pi)\geq (1-\varepsilon)\omega(T),$$ where $T_{\pi}=E_{\pi} T E_{\pi} \in Lip_0(\OY_{\pi})$. We can identify  $\OY_{\pi}$ with $l_1^n(\OX)$ for some $n\in \mathbb{N}$. Then $T_\pi\in Lip_0(\OY_{\pi})$ is given by $T_\pi(x)=(Sx_1, Sx_2, \cdots, Sx_n)$ for each $x=(x_1, x_2, \cdots, x_n)$ in $l_1^n(\OX)$. By Lemma \ref{finite range}, we have $\omega(T_\pi)=\omega(S)$ and thus $$\omega(S)\geq (1-\varepsilon)\omega(T)\geq (1-\varepsilon) n_L(L_1(\mu, \OX)).$$ This completes the proof. {\hfill $\Box$\medskip}
\medskip

To generalize the fact $n_L(l_\infty(\OX))= n_L(\OX)$, we require the following two results of \cite[Lemmas 1, 2]{MV}. For the convenience of the readers, we include them here.

\begin{lem}\label{preserve}
Let $(\Omega,\Sigma,\nu)$ be a $\sigma$-finite measure space, and let $\OX$ be a Banach space.
If $f\in L_\infty(\nu, \OX)$ with $\|f(t)\|>\lambda$ a.e.,  then there exists $B\in \Sigma$ with $0<\nu(B)<\infty$ such that
$$\|\frac{1}{\nu(B)}\int_B f(t) d\nu(t)\|>\lambda.$$
\end{lem}

\begin{lem}\label{dense}
Let $f\in L_\infty(\nu, \OX)$, $C\in \Sigma$ with positive measure and $\varepsilon>0$. Then there exist $x \in \OX$ and $A \subset C$ with
$0<\nu(A)<\infty$ such that $\|x\|=\|f\chi_c\|$ and $\|(f-x)\chi_A\|<\varepsilon$. Accordingly, the set
$$\{x\chi_A+f\chi_{\Omega/A}: x\in S_\OX, f\in B_{L_\infty(\nu, \OX)}, A\in \Sigma \ with \ 0<\nu(A)<\infty\}$$
is dense in $S_{L_\infty}(\nu, \OX)$.
\end{lem}

\begin{thm}\label{finite measure}
Let $(\Omega,\Sigma,\nu)$ be a $\sigma$-finite measure space, and let $\OX$ be a Banach space. Then
$$n_L(L_\infty(\nu, \OX))=n_L(\OX).$$
\end{thm}

\begin{prf}
We first show that $n_L(L_\infty(\nu, \OX))\geq n_L(\OX).$ Fix $T\in Lip_0(L_\infty(\nu, \OX))$ with $\|T\|_L=1$. We need to prove that $\omega(T)\geq n_L(\OX)$. Given $\varepsilon>0$, there are $f,g\in L_\infty(\nu, \OX)$
and a set $C\in \Sigma$ with $0<\nu(C)<\infty$ such that $$\|Tf(t)-Tg(t)\|>(1-\varepsilon/2)\|f-g\|, \ \forall\  t\in C.$$
Denote $f_0=\frac{f-g}{\|f-g\|}$ and set $$\mathcal{A}(f_0)=\{x\chi_A+f_0\chi_{\Omega/A} \in S_{L_\infty(\nu, \OX)}:  x\in S_\OX, A\in \Sigma, A\subset C,  \nu(A)>0 \}.$$ By the proof of \cite[Theorem 3]{MV}, we have $f_0 \in \overline{J(\mathcal{A}(f_0))}$ (Indeed, by Lemma \ref{dense}, there exist $y_0 \in B_\OX$ and $A \subset C$ with $\nu(A) > 0$ such that $\|(f_0-y_0)\chi_A\| <\varepsilon.$ Now, write $y_0 = \lambda x_1 +(1-\lambda)x_2$ with $0\leq \lambda \leq 1, x_1, x_2 \in S_\OX$, and consider the functions $$f_j = x_j \chi_A + f_0 \chi_{\Omega/A} \in \mathcal{A}(f_0) \quad (j =1, 2),$$ which clearly satisfy $\|f_0 -(\lambda f_1 +(1-\lambda)f_2)\| <\varepsilon$). From Lemma \ref{join}, we can find $h\in L_\infty(\nu, \OX)$, $z\in S_\OX$ and $A_1\in \Sigma$ with $A_1\subset C$  and $0<\nu(A_1)<\infty$ such that $$ h-f= \frac{\|f-g\|}{2}(z\chi_{A_1}+f_0\chi_{\Omega/A_1})\in \frac{\|f-g\|}{2}\mathcal{A}(f_0) \quad \mbox{and} \quad \|h-g\|\leq (1+\varepsilon)\frac{\|f-g\|}{2}.$$ Hence, we have \begin{eqnarray}\|T h(t)-T f(t)\|\geq \|T f(t)-T g(t)\|-\|h-g\|\geq (1-2\varepsilon)\|h-f\|, \ \forall\  t\in C.\end{eqnarray} By Lemma \ref{dense}, there exist $x_0, y_0\in B_\OX$ and $A_0\subset A_1$ with $0<\nu(A_0)<\infty$ such that $$\|(h-x_0)\chi_{A_0}\|<\frac{\varepsilon}{2}\|h-f\|  \quad \mbox{and} \quad  \|(f-y_0)\chi_{A_0}\|<\frac{\varepsilon}{2}\|h-f\|.$$  Note that for each $t\in A_0\subset A_1$,  $\|h-f\|=\|h(t)-f(t)\|$. Therefore, for all $t\in A_0$, $$\|x_0-y_0\|\geq \|h(t)-f(t)\|-\|h(t)-x_0\|-\|f(t)-y_0\|>\|h-f\|(1-\varepsilon)$$ and
$$\|x_0-y_0\|\leq \|h(t)-f(t)\|+\|h(t)-x_0\|+\|f(t)-y_0\|<\|h-f\|(1+\varepsilon).$$ This means that $$\|h-f\|(1-\varepsilon)<\|x_0-y_0\|<\|h-f\|(1+\varepsilon).$$
By Lemma \ref{extension}, there exists a $\frac{1}{1-\varepsilon}$-Lipschitz operator $F\in Lip(\OX, L_\infty(\nu, \OX))$ such that $$F(x_0)=h\ \ \mbox{and} \ \ \ F(y_0)=f.$$ For each $A\in \Sigma$ with $0<\nu(A)<\infty$, we define a contractive linear operator $\Theta_A: L_\infty(\nu, \OX)\rightarrow \OX$ by $$\Theta_A(f)=\frac{1}{\nu(A)}\int_A f d\nu, \quad  \forall f\in L_\infty(\nu, \OX).$$
Now by (9) and Lemma \ref{preserve} we select $B\subset A_0\subset C$ with $0<\nu(B)<\infty$ satisfying $$\|\Theta_B (Th-Tf)\|\geq(1-2\varepsilon)\|h-f\|.$$
We denote $$\Phi(x)=x\chi_{A_0}+F(x)\chi_{\Omega/A_0} \in L_\infty(\nu, \OX),$$ and consider the operator $S \in Lip_0(\OX)$ given by $$S(x)=\Theta_B [T(\Phi(x))-T(\Phi(0))].$$
Since \begin{eqnarray*}
\|S(x_0)-S(y_0)\|&=&\|\Theta_B [T(x_0\chi_{A_0}+h\chi_{\Omega/A_0})]- \Theta_B[T(y_0\chi_{A_0}+f\chi_{\Omega/A_0})] \|
\\ &\geq&\|\Theta_B(T h-T f)\|-\|(h-x_0)\chi_{A_0}\|- \|(f-y_0)\chi_{A_0}\|
\\ &\geq&(1-3\varepsilon)\|h-f\|\geq\frac{(1-3\varepsilon)}{1+\varepsilon}\|x_0-y_0\|>(1-4\varepsilon)\|x_0-y_0\|,
\end{eqnarray*} this implies that $$\|S\|_L\geq 1-4\varepsilon.$$ For each $x, y\in \OX$ and $(x-y)^*\in D(x-y)$, we have
$$(x-y)^*\circ \Theta_B (\Phi(x)-\Phi(y))=\|x-y\|^2$$
and $$(x-y)^*\circ \Theta_B(T\Phi(x)-T\Phi(y))=(x-y)^*(Sx-Sy).$$
It follows that $$\omega(T)\geq \omega(S)\geq n_L(\OX)(1-4\varepsilon)$$ and so $n_L(L_\infty(\nu, \OX))\geq n_L(\OX)$.
\medskip

To prove the reverse inequality, for every $S \in Lip_0(\OX)$ with $\|S\|_L=1$, we define $T\in Lip_0(L_\infty(\nu, \OX))$ by
$(Tf)(t)=S(f(t))$ for each $t\in \Omega, \ f\in L_\infty(\nu, \OX).$ Then $\|T\|_L=1$ and $\omega(T)\geq n_L(L_\infty(\nu, \OX))$.
Lemma \ref{dense} combined with Corollary \ref{choose} produces $x, y\in X$ with $x\neq y$, $(x-y)^*\in D(x-y)$ and $f, g\in L_\infty(\mu,X)$ in the form $f=x\chi_A+ f_0\chi_{\Omega\setminus A}$, $g=y\chi_A+ g_0\chi_{\Omega\setminus A}$ for some $A\in \Sigma$ with $0<\nu(A)<\infty$ and some $f_0,g_0\in L_\infty(\mu,X)$ such that
$$\omega(T)-\varepsilon<\frac{|(x-y)^*\circ \Theta_A(Tf-Tg)|}{\|f-g\|^2}\leq \frac{|(x-y)^*(Sx-Sy)}{\|x-y\|^2}.$$
Therefore, $$n_L(L_\infty(\nu, \OX))-\varepsilon\leq \omega(T)-\varepsilon\leq \omega(S),$$ and thus
$n_L(L_\infty(\nu, \OX))\leq n_L(\OX).$
\end{prf}

The following result is an obvious consequence of Theorems \ref{compact}, \ref{positive measure} and \ref{finite measure}.
\begin{cor}
Let $K$ be a compact Hausdorff space, $\mu$ be a positive measure and $\nu$ be a $\sigma$-finite measure.
Then, in the real or complex case, one has
$$n_L(C(K))=n_L(L_1(\mu))=n_L(L_\infty(\nu))=1.$$
\end{cor}

Throughout this paper, we can see that the conclusion that the Lipschitz numerical index is
equal to the numerical index is true for large classes of Banach spaces. Then it is interesting to ask the following question:
\begin{problem}
Is there a Banach space $\OX$ such that $n_L(\OX)<n(\OX)$?
\end{problem}

\subsection*{Acknowledgments}The authors would like to thank Guanggui Ding for fruitful conversations
concerning the matter of this paper.

\end{document}